\documentclass[11pt]{article}

\usepackage{amsmath}
\usepackage{amsfonts}
\usepackage{theorem}

\evensidemargin \oddsidemargin
\def\R{{\mathbb R}}
\def\N{{\mathbb N}}
\def\Z{{\mathbb Z}}
\def\P{{\mathbb P}}

\def\Proof{{\smallskip\noindent{\em Proof. }}}
\def\endProof{{\hfill$\bullet$\medskip\noindent}}

\theoremstyle{plain}
\newtheorem{theorem}{Theorem} [section]

\newtheorem{lemma}[theorem]{Lemma}

\newtheorem{corollary}[theorem]{Corollary}

\newtheorem{proposition}[theorem]{Proposition}
{\theoremstyle{plain}
  \theorembodyfont{\normalfont\rmfamily}

\newtheorem{remark}[theorem]{Remark}

            }

\begin{document}

\title{Space-time decay of Navier--Stokes\\
       flows invariant under rotations}

\date{March, 15 2003}

\author{Lorenzo Brandolese\\
	\\
        Centre de Math\'ematiques et de Leurs applications\\
	ENS de Cachan\\
	61, Avenue du Pr\'esident Wilson\\
	94235 Cachan Cedex, FRANCE\\
	e-mail: brandole@cmla.ens-cachan.fr      }

\maketitle

\abstract{
We show that the solutions to the non-stationary Navier--Stokes equations in $\R^d$ $(d=2,3)$
which are left invariant under the action of discrete subgroups of the orthogonal group $O(d)$
decay much faster as $|x|\to\infty$ or $t\to\infty$ than in generic case and we compute, for each subgroup,
the precise decay rates in space-time of the velocity field.
%
}

\section{Introduction and main results}
This paper is devoted to the study of the asymptotic behavior of viscous flows
of incompressible fluids filling the whole space $\R^d$ ($d\ge2$) and not submitted to the action
of external forces.
These flows are governed by the Navier--Stokes equations, which we may write in the following form 
$$
\left\{
\begin{array}{l}
\partial_t u -\Delta u + \P\nabla\cdot(u\otimes u)  = 0 \\
u(x,0)=a(x)\\
\nabla\cdot u=0.
\end{array}
\right.
\eqno\hbox{(NS)}
$$
Here $u(\cdot,t)\,:\,\R^d\to\R^d$ $(d\ge2)$ denotes the velocity field and
$\P$ is the Leray--Hopf projector onto the soleinodal vectors field, defined by
$\P f= f - \nabla\, \Delta^{-1}(\nabla\cdot f)$, with $f=(f_1,\ldots,f_d)$.

It is now well known (see {\em e.g.\/} \cite{DS94}, \cite{Bra01}, \cite{Sch91}, \cite{GW2})
that {\em generic\/} solutions $u$ to (NS) decay at infinity at considerably slow rates
in space-time.
Indeed, even if the data have the form $a(x)=\epsilon \phi(x)$, 
where $\epsilon>0$ is a small constant, the components of $\phi$ belong to the Schwartz class
and have vanishing moments, then the corresponding strong solution $u$ to (NS)
satisfies
$|u(x,t)|\le C(1+|x|)^{-(d+1)}$ and $|u(x,t)|\le C(1+t)^{-(d+1)/2}$ for all $x\in\R^d$ and $t>0$;
but such decay rates, in general, are optimal
(we refer to \cite{Mi00}, \cite{Mi02II} \cite{AGSS}, \cite{HX}
for a proof of these bounds under different assumptions).
Furthermore, very few examples of solutions which decay faster
are known so far: we should mention here the classical example of a two dimensional flow with
{\em radial vorticity\/} (\cite{Sch91}, \cite{FM}, see also \cite{SchSS}) and the examples of flows
in $\R^d$ ($d\ge2$) constructed in \cite{Bra01}, \cite{Bra02} imposing some special symmetries on
the initial data.

The purpose of this paper is to provide a systematic study of the connection between symmetry and
space-time decay of viscous flows in dimension two and three.
Our starting point is the observation that the Navier--Stokes equations are invariant under
the transformations of the orthogonal group $O(d)$:
if $u(x,t)$ is a solution to the Navier--Stokes equations in
$\R^d$, and $P\in O(d)$ is an orthogonal matrix, then $\tilde u(x,t)=P^Tu(P\,x,t)$
is a Navier--Stokes flow as well (here, $P^T$ is the transposed matrix).
Roughly, it follows that if the initial data commute with $P\in O(d)$,
then the velocity field will satisfy
\begin{equation}
\label{inv-u}
P\,u(x,t)=u(P\,x,t),
\end{equation}
whenever the solution to (NS) is defined (even if only in a weak sense).
 
If $G$ is any subgroup of $O(d)$, then a natural problem is that of computing 
the space-time decay rates of solutions that are invariant under all the transformations of~$G$.
In this paper we will consider only {\em discrete\/} subgroups of the orthogonal group, the reason being the
following: in the two dimensional case, solutions that are invariant under the continuous
subgroup $SO(2)$ do exist, but boil down to flows with {\em radial vorticity\/}.
These flows are ``trivial''
in the sense that the non-linear term $\P\nabla\cdot(u\otimes u)$ in (NS)
identically vanishes.
On the other hand, in the three dimensional case, one easily sees via the Fourier
transform that flows $u(x,t)\not\equiv0$ which are invariant under the whole group $SO(3)$ do not exist.

As an immediate consequence of our results in the $d=2$ case
we will prove in section~\ref{section2D} the following theorem.

\begin{theorem}
\label{theorem1}
Let $a(x)$ be a soleinoidal vector field in ${\cal S}(\R^2)$ (the Schwartz class).
If $a$ is invariant under the cyclic group $C_n$ of order $n$, then the global strong solution
$u(x,t)$ such that $u(0)=a$ satisfies $u(x,t)=O(|x|^{-(n+1)})$ for all $t\ge0$.
If, in addition, $a$ is invariant under the dihedral group $D_n$ of order $2n$, then this decay is uniform 
in time and, moreover, $||u(t)||_p\le C(1+t)^{-(n+1)/2+1/p}$ (with $2<p\le\infty$).
\end{theorem}

As we shall see later on, the flows obtained in Theorem~\ref{theorem1} 
do not have radial vorticity and hence they do not boil down to ``trivial solutions''
(or to solutions of the homogeneous heat equation $\partial_t u=\Delta u$).
At best of our knowledge, no other examples of highly localized flows in $\R^2$ were known so far.

\bigskip
In the three dimensional case,
the problem of the existence of (trivial or non-trivial) rapidly decreasing solutions 
$u=(u_1,u_2,u_3)\not\equiv0$ as $|x|\to\infty$ to
Navier--Stokes equations was raised in \cite{DS96} and it is still open.

However, non-trivial and localized divergence-free vector fields in $\R^3$ $a(x)$,
which are invariant under discrete subgroups of $O(3)$ can be easily constructed, and we may expect
that such fields should lead to solutions with fast decay at infinity.
If $G$ is one of these subgroups, then we know that $G$ is either
\begin{itemize}
\item a subgroup of the complete symmetry group of a regular polyhedron, or
\item a subgroup of the complete symmetry group of a prism (and hence isomorphic to a cyclic or a dihedral
group).
\end{itemize}
As we will see in section~\ref{optimality},
solutions which are invariant under the complete symmetry group of a prism
(which is not a cube) in general do not decay faster than $|x|^{-4}$
(the same remark applies for continuous subgroups of $SO(3)$ such as the complete direct symmetry
group of the cylinder).
On the other hand, flows with polyhedral symmetry decay much faster.
As a byproduct of our constructions, we shall be able to provide examples of solutions $u(x,t)$ 
decaying at infinity as $|x|^{-8}$ and $t^{-4}$, thus improving the results of \cite{Bra01},
\cite{Mi02} and \cite{Bra02}.
The most interesting cases are described in the theorem below (see section~\ref{section3D} for
the exaustive study of the asymptotic behavior of all the other finite groups of isometries in $\R^3$).

Let us denote by $L^\infty_\gamma(\R^d)$ $(\gamma\ge0)$ 
the space of all measurable functions (or vector fields) $f$,
defined on $\R^d$, and such that $(1+|x|)^\gamma|f(x)|\in L^\infty(\R^d)$.
For any positive $T$, $0<T\le\infty$, we denote by $C([0,T],L^\infty_\gamma(\R^d))$ the
space of continuous and bounded $L^\infty_\gamma(\R^d)$-valued functions, the continuity
at $t=0$ being understood in the distibutional sense.
Then we have the following:

\begin{theorem}
\label{theorem2}
Let $a=(a_1,a_2,a_3)$ a divergence-free and rapidly decreasing vector field in $\R^3$:
$(1+|x|)^k a\in L^\infty(\R^3)$, for all $k=0,1,\ldots$.
Then we know (\cite{Mi00}, \cite{BrM}) that there exists $T$  $(0<T\le\infty)$ and a unique strong
solution $u$ to the Navier-Stokes
equations in $\R^3$, such that $u(0)=a$ and $u\in C([0,T],L^\infty_4(\R^3))$.
\begin{enumerate}
\item
If $a(x)$ is invariant under the complete symmetry group of the tetrahedron, then
$u\in C([0,T],L^\infty_5(\R^3))$.
\item
If $a(x)$ is invariant under the complete symmetry group of the cube (or of the octahedron), then
$u\in C([0,T],L^\infty_6(\R^3))$.
\item
If $a(x)$ is invariant under the complete symmetry group of the dodecahedron (or of the icosahedron),
then $u\in C([0,T],L^\infty_8(\R^3))$.
\end{enumerate}
Furthermore,
if we know that $u(x,t)$ is global ($T=\infty$), then $||u(t)||_p$ decays, 
respectively, at least as fast as $t^{-5/2+3/(2p)}$, $t^{-3+3/(2p)}$ and $t^{-4+3/(2p)}$ 
 as $t\to\infty$ $({3\over 2}< p\le\infty)$.
\end{theorem}

\begin{remark}
The result of Theorem~\ref{theorem2} is sharp.
Optimality of the above space decay rates should be understood in the following sense:
if $G$ is one of the previous three groups and $\gamma=5$, $6$ or $8$ (respectively),
then there exists a solution $u(x,t)$ to (NS) which is invariant under $G$ and  localized at $t=0$,
but which does not decay faster than $|x|^{-\gamma}$, uniformely in any time interval
$[0,\epsilon]$ $(\epsilon>0)$.
For each group $G$ we shall provide examples of such flows.
\end{remark}

We shall see in section~\ref{section3D} that, because of the symmetries imposed on the initial data,
the velocity field has vanishing moments 
$\int x^\alpha u(x,t)\,dx$ up the the order~$1$,
$2$ and~$4$, for all $t\in[0,T]$, respectively in the case~1,~2 and 3 of Theorem~\ref{theorem1}.
In particular, the fact that $a(x)$ has these cancelations allows us to see
that the estimates in space-time obtained for $u(x,t)$ hold true also for the linear
evolution $e^{t\Delta}a(x)$ (here $e^{t\Delta}$ denotes the heat semigroup).

Let us point out that the existence of a {\em global\/} strong solution is usually ensured
by some smallness assumption on the initial data: a common suitable assumption is {\it e.g.\/} that
$||a||_3$ is small enough (see \cite{Ka84}).
However, for the flows treated in Theorem~\ref{theorem2}, the three equations contained in the first of (NS)
reduce to a simpler {\em single scalar equation\/} on the first component $u_1(x,t)$.
Thus, it would be an interesting problem to study the global solvability of those ``symmetric'' solutions
in the case of ``large'' initial data.

\medskip
There is an extensive literature on the asymptotic behavior of the Navier--Stokes equations 
(see {\it e.g.} \cite{HX}, \cite{GW1}, \cite{GW2}, \cite{Mi00}, \cite{Sch91}, \cite{Wie87}
and the references therein contained),
but not so much has been written on symmetry of viscous flows.
See, however, \cite{Kid85I}, \cite{Kid85II} for applications
of symmetries to the numerical simulation of turbulence and \cite{FP94}, \cite{Pop95} for the construction
of ansatzes to (NS). 
The connection between symmetry and space-time decay have been first noticed in \cite{Bra01} and
subsequentely studied in \cite{Bra02}, \cite{Mi02}.
The symmetries which are considered in these papers are only those corresponding to a 
subgroup of the group of the symmetries of the cube.
Hence, the main results of \cite{Bra01}, \cite{Mi02}, \cite{Bra02} 
are contained in the present paper as a particular case.

It is worth observing that recently Th.~Gallay and C.\,E.~Wayne were able to prove
the existence of flows with a fixed, but arbitrarily large, {\it time\/} decay rate
(see \cite{GW1} and \cite{GW2}, respectively for $d=2,3$).
Indeed, using the vorticity formulation of the Navier--Stokes equations, they showed in
\cite{GW1} and \cite{GW2} that the class of solutions which decay faster than a given rate as $t\to\infty$
lies on an invariant manifold of finite codimension, in a suitable functional space.
Their method, which is a combination of the spectral decomposition of the
Fokker--Planck operator and the theory of dynamical systems, would
be effective in any space dimension.
However, this approach yields no {\em explicit\/} examples of initial data leading to such 
solutions with fast decay.

\medskip
The rest of this paper is organized as follows.
In section~\ref{rappels} we briefly recall the vorticity formulation of (NS) and 
a general result of the author on the space decay
of solutions to the Navier--Stokes equations that we will use throughout this paper. 
As an application of this result to the two-dimensional case, in section~\ref{section2D}
we will prove Theorem~\ref{theorem1} in a slightly more general form.
In section~\ref{section3D} we start recalling the complete list of the discrete subgroups of $O(3)$
and we subsequentely compute the space decay rates of flows invariant under the action of all these groups.
There we will also discuss the closely related problem of the cancellations of the vorticity of such flows.
In section~\ref{optimality} we will show  by means of some examples the optimality
of the decay rates that we obtain.

\section{Decay of the velocity field and the vorticity}
\label{rappels}

Throughout this paper we shall assume that the initial datum $a$ is a rapidly decreasing function in
$\R^d$ $(d\ge2)$.
This requirement is not essential (the optimal assumptions should be expressed in terms of Besov
and weak-Hardy spaces, as in \cite{Mi00}, \cite{Mi02II} and \cite{Bra02})
but it considerably simplifies the presentations of our main results.
Then we know that a {\em necessary condition\/} on the data, in order to avoid that the velocity field
instantaneously ``spreads out'', is that the components of $a$ are orthogonal with respect to the
$L^2$ inner product (see \cite{DS94}):
\begin{equation}
\label{orth}
\int (a_ha_k)(x)\,dx= c\delta_{h,k}, \qquad (h,k=1,\ldots,d)
\end{equation}
($\delta_{h,k}=1$ if $h=k$ and $\delta_{h,k}=0$ if $h\not=k$). 
More in general, using the Fourier transform it is not difficult to see that if we want that the solution
$u(t)$ remains rapidly decreasing as $|x|\to\infty$, at least in a small time interval $[0,T]$ $(T>0)$,
then for all $m=0,1,\ldots$ the homogeneous polynomial $P_m(a)$, defined by
\begin{equation}
\label{Pma}
P_m(a)(\xi)\equiv \sum_{h,k=1}^d\sum_{|\alpha|=m} \biggl({1\over \alpha!}\int x^\alpha(a_ha_k)(x)\,dx \biggr)
  \xi^\alpha\xi_h\xi_k,
\end{equation}
(here $\xi=(\xi_1,\ldots,\xi_d)\in\R^d$ and we adopted the usual notations for the multi-index 
$\alpha=(\alpha_1,\ldots,\alpha_d)\in\N^d$)
{\em must be divisible\/} by $\xi_1^2+\cdots+\xi_d^2$.

Conversely, if we assume that this last condition is satisfied at the beginning of the evolution
and that {\it it remains true\/} during a time interval $[0,T]$, then the corresponding strong 
solution will be rapidly decreasing as $|x|\to\infty$ for all $t$ in such interval.
More precisely, let us recall the following result of \cite{BrM}, which will be our main tool
for our study of the spatial localization. 

\begin{proposition}
\label{proposition1}
Let $M$ be a fixed non-negative integer and $a(x)$ a divergence-free and rapidly decreasing vector field
in $\R^d$.
Let also $u(x,t)$ be the unique strong solution to (NS), defined in some time interval $[0,T]$ 
$(0<T\le\infty)$ such that $u(0)=a$ and $u\in C([0,T],L^\infty_{d+1}(\R^d))$.

If the polynomials $P_m(u(t))(\xi)$, defined as in \eqref{Pma}, are divisible by $\xi_1^2+\cdots+\xi_d^2$
for all $t\in[0,T]$ and $m=0,1,\ldots,M$, then the spatial decay of $u$ is improved by
$ u\in C([0,T'],L^\infty_{d+2+M}(\R^d))$,
for any $T'\in \R^+$ $(0\le T'\le T)$.
Furthermore, if $T=+\infty$ and the moments of $a$ vanish up to the order $1+M$, then
\begin{equation}
\label{decay-u-infinity}
 u\in C([0,+\infty[,L^\infty_{d+2+M}(\R^d)),
\end{equation}
\end{proposition}

\medskip
The assumpion on the moments of $a$ ensures that, if the solution is globally defined, then
the norm $||u(t)||_{L^\infty_{d+2+M}}\equiv \sup_x (1+|x|)^{d+2+M}|u(x,t)|$ does not blow up
as $t\to\infty$.
Condition~\eqref{decay-u-infinity} does not ensures that the solution decay fast as $t\to\infty$.
It can be shown, however, that if the moments of $a$ vanish up to the order $1+M$
{\em and if the following identities hold true\/}:
\begin{eqnarray}
\label{iden-time-decay}
\lefteqn{ \sum_{h,k=1}^d\sum_{|\alpha|=m} \biggl({1\over \alpha!}\int x^\alpha(u_hu_k)(x,t)\,dx \biggr)
  \xi^\alpha\xi_j\xi_h\xi_k \equiv} \\
&& \sum_{h=1}^d\sum_{|\alpha|=m} 
       \biggl({1\over \alpha!}\int x^\alpha(u_hu_j)(x,t)\,dx  \biggr)\xi^\alpha\xi_h (\xi_1^2+\cdots+\xi_d^2),
       \qquad(j=1,\ldots,d) \nonumber
\end{eqnarray}
for all $\xi\in\R^d$, $t\ge0$ and $m=0,1,\ldots,M$,
then
\begin{equation}
\label{xt-decay-u}
\sup_{x\in\R^d,t\ge0}  (1+|x|)^{\gamma}(1+t)^{(d+2+M-\gamma)/2}|u(x,t)|<\infty \qquad(0\le\gamma\le d+2+M).
\end{equation}
We refer to \cite{Bra02} for a proof of this claim and examples of flows satisfying~\eqref{iden-time-decay}
for $m=0,1$.
See also \cite{FM}, \cite{MS}, \cite{GW2} for related results.
Condition~\eqref{iden-time-decay}, however, is difficult to check for large $m$.
To construct examples of flows with fast decay in space-time we shall rather make use
of the vorticity formulation of the Navier--Stokes equations.
This allows us to give a much more natural sufficient condition which ensures~\eqref{xt-decay-u}.

\medskip
From now on we shall work only in two or three space dimension.
We recall that the  vorticity is defined by
$$
\omega=\partial_1 u_2-\partial_2 u_1 \qquad \hbox{($d=2$)}
$$
or,
\begin{equation}
\Omega=\nabla\times u \nonumber
           =(\partial_2 u_3-\partial_3 u_2,\partial_3 u_1-\partial_1 u_3,
	      \partial_1 u_2-\partial_2 u_1) \qquad(d=3).
\end{equation}
Note that the vorticity is a scalar function when $d=2$ and a soleinoidal vector field if $d=3$.
Then the vorticity verifies the integro-differential equations
\begin{equation}
\label{eq-omega}
\partial_t\omega+(u\cdot\nabla)\omega =\Delta \omega 
         \qquad(d=2)
\end{equation}
or,
\begin{equation}
\label{eq-Omega}
\partial_t\Omega+(u\cdot\nabla)\Omega-(\Omega\cdot\nabla) u =\Delta\Omega, \qquad\nabla\cdot\Omega=0
	\qquad(d=3).
\end{equation}
Here, the velocity field $u$ has to be expressed in terms of its vorticity via the Biot--Savart laws:
\begin{align}
&u(x,t)={1\over 2\pi} \int {(x-y)^\perp\over |x-y|^2} \omega(y,t)\,dy,  \qquad\quad(d=2), \label{BS2}\\
&u(x,t)=-{1\over 4\pi} \int { x-y\over |x-y|^{3}}\times \Omega(y,t)\,dy \qquad(d=3), \label{BS3}
\end{align}
where we denoted $(x_1,x_2)^\perp=(-x_2,x_1)$ in~\eqref{BS2}.

We collect in the following proposition several known facts on the vorticity equation:

\begin{proposition}
\label{proposition2}
\begin{enumerate}
\item
Let $n$ be a positive integer, $\omega_0$ a rapidly decreasing function in $\R^2$
with vanishing moments up to the order $n-1$ and $\omega(x,t)$ the  unique global strong solution
of~\eqref{eq-omega}-\eqref{BS2}, such that $\omega(0)=\omega_0$.
If we know that the moments of $\omega(t)$ vanish up to the order $n-1$, for all $t\ge0$, then
we have
\begin{equation}
\label{profiles-omega}
\sup_{x,t}\,(1+|x|)^\gamma(1+t)^{(2+n-\gamma)/2}|\omega(x,t)|<\infty
\end{equation}
for all $\gamma\ge0$ and 
\begin{equation}
\label{Lp-decay-omega}
||\omega(t)||_p\le C (1+t)^{-(2+n)/2+1/p} \qquad(1\le p\le\infty).
\end{equation}

\item
Let $\Omega_0$ be a rapidly decreasing and divergence-free vector field in $\R^3$,
with vanishing moments up to the order $n-1$.
If $\sup_x |x|^2|\Omega_0(x)|$ is small, then there exists a unique 
strong solution $\Omega(x,t)$ of~\eqref{eq-Omega}-\eqref{BS3}, such that $\Omega(0)=\Omega_0$
and $\Omega\in C([0,\infty[,L^\infty_2(\R^3))$.
If we know that the moments of $\Omega(t)$ vanish up to the order $n-1$, for all $t\ge0$, then
we have
\begin{equation}
\label{profiles-Omega}
\sup_{x,t}\,(1+|x|)^\gamma(1+t)^{(3+n-\gamma)/2}|\Omega(x,t)|<\infty
\end{equation}
for all $\gamma\ge0$ and
\begin{equation}
\label{Lp-decay-Omega}
||\Omega(t)||_p\le C (1+t)^{-(3+n)/2+3/(2p)} \qquad(1\le p\le\infty).
\end{equation}
\end{enumerate}
\end{proposition}

\begin{remark}
We refer to \cite{BA94} and  \cite{GM} for the study of the well-posedness of
the Cauchy problem for equations~\eqref{eq-omega} and~\eqref{eq-Omega}.
In particular, it is well known that \eqref{eq-omega}-\eqref{BS2} can be uniquely solved {\it e.g.\/} in
$C([0,\infty),L^1(\R^2))\cap C(]0,\infty),L^\infty(\R^2))$.
In the three dimensional case (and in the case of small initial data),
the fact that \eqref{eq-omega}-\eqref{BS2}
can be uniquely solved in $C([0,\infty[,L^\infty_2(\R^3))$ is easily seen, see \cite{Bra03}.

Note that the decay profiles of $\omega$ and $\Omega$ 
are the same which can  be obtained for the solutions
of the homogeneous heat equations $e^{t\Delta}\omega_0(x)$ and $e^{t\Delta}\Omega_0(x)$, respectively.
We refer to \cite{Bra03} for a proof of~\eqref{profiles-Omega} (the proof of~\eqref{profiles-omega} is identical).
The decay of the vorticity in the $L^p$-norm are formally a consequence of~\eqref{profiles-omega} and~\eqref{profiles-Omega},
respectively if $d=2$ or $3$.
Estimates~\eqref{Lp-decay-omega} and~\eqref{Lp-decay-Omega}, however, can be proved with straightforward adaptations
of the arguments of \cite{FM}-\cite{Bra03}, or \cite{GW1}-\cite{GW2} 
(see also \cite{Car96}).
\end{remark}

We would like to stress the fact that, 
since the moments of the vorticity are {\it not invariant\/} during the time evolution
(excepted for the integral and the first order moments), profiles~\eqref{profiles-omega} and~\eqref{profiles-Omega}
will hold true only for $n=0,1,2$ even if the vorticity has many cancelations at time $t=0$.
This reflects the fact that, in general, the velocity field does not decay faster
than $|x|^{-(d+1)}$.
We will show, however, that for the special flows described in Theorem~\ref{theorem1} and Theorem~\ref{theorem2}
the vorticity has a large number of vanishing moments for all time
(in particular, we shall be able to prove that the assumptions of the second part of Proposition~\ref{proposition2}
are non-vacuous for $n=0,\ldots,6$).

We will finish this section stating another simple result
which allows us to deduce time decay estimates for the velocity field from
\eqref{profiles-omega}-\eqref{profiles-Omega}.

\begin{lemma}
\label{lemma1}
Let $\omega(x,t)$ ($d=2$), or  $\Omega(x,t)$ $(d=3)$, as in Proposition~\ref{proposition2},
and let $u(x,t)$ be the corresponding velocity field obtained via the Biot--Savart law
\eqref{BS2} (resp.~\eqref{BS3}).
Then we have,
\begin{equation}
\label{Lp-decay-u-2d}
||u(t)||_p \le C(1+t)^{-(n+1)/2+1/p}, \qquad 2<p\le \infty \qquad(d=2),
\end{equation}
or
\begin{equation}
\label{Lp-decay-u-3d}
||u(t)||_p \le C(1+t)^{-(n+2)/2+3/(2p)}, \qquad 3/2 <p\le \infty \qquad(d=3)
\end{equation}
and the spatial moments of $u(x,t)$ exist and vanish up to the order $n-2$.
\end{lemma}

\Proof
Note that the Biot--Savart kernels $x^\perp/|x|^{2}$ $(d=2)$ and $x/|x|^3$ $(d=3)$ belong respectively
to the the weak-Lebesgue spaces $L^{2,\infty}(\R^2)$ and $L^{3/2,\infty}(\R^3)$.
Bounds \eqref{Lp-decay-u-2d} and \eqref{Lp-decay-u-3d} then immediately follow from the Biot--Savart laws \eqref{BS2}
and \eqref{BS3},
the corresponding bounds for the vorticity and elementary results on convolution and interpolation of Lorentz spaces
(see \cite{BL}).
The condition on the moments of $u$ is easily seen via the Fourier transform and the Taylor formula
(see {\it e.g.\/} \cite{Bra03} and \cite{GW2} for this type of calculation).
\endProof

\section{Space-time decay of two-dimensional flows}
\label{section2D}

Let us recall that all finite subgroups of the orthogonal group $O(2)$, are of two types:
cyclic groups (which are indeed subgroups of the special orthogonal group $SO(2)$ of proper rotations),
and dihedral groups.
We shall denote by $C_n$ the cyclic group of order $n$ and by $D_n$ the dihedral group of order $2n$.
This group contains $C_n$ and its presentation is given by
two generators $R$ and $\tau$, together with the relations $R^n={\bf 1}$, $\tau^2={\bf 1}$ and $\tau R=R^{-1}\tau$
($R$ corresponds to a rotation of $2\pi/n$ around the origin 
and $\tau$ to a reflection with respect to a straght line passing through the origin).

Divergence-free vector fields, which are rapidly decreasing as $|x|\to\infty$ and
which are left invariant under the actions of $C_n$ or $D_n$ are easily constructed by means of the
vorticity.
Indeed, in general, if $P\in O(2)$, $u(Px)=Pu(x)$ for all $x\in\R^2$ and $\omega=\partial_1u_2-\partial_2u_1$,
then $\omega(x)=\det(P)\omega(Px)$ (in the distributional sense).
Conversely, if $\omega(x)=\det(P)\omega(Px)$ for all $x$ and $u$ is given by the Biot--Savart law,
then $u(Px)=Pu(x)$, whenever the singuar integral \eqref{BS2} makes sense.

\medskip
Assume now that the initial datum $a=u(0)$ is rapidly decreasing as $|x|\to\infty$ and
that it is {\em invariant under the cyclic group of order $n$\/} $(n\ge3)$.
Since we know that strong solutions of the two dimensional Navier--Stokes equations are globally defined,
we have $u(Rx,t)=Ru(x,t)$, for all $x\in\R^2$ and $t\ge0$.
We now apply to $u$ Proposition~\ref{proposition1}:
let us show that the polynomial
\begin{equation}
\label{Pmu}
P_m(u)(\xi)\equiv \sum_{h,k=1}^2\sum_{|\alpha|=m} \biggl({1\over \alpha!}\int x^\alpha(u_hu_k)(x,t)\,dx \biggr)
  \xi^\alpha\xi_h\xi_k,
\end{equation}
is divisible by $\xi_1^2+\xi_2^2$, for all $t\ge0$ and the first values of $m$.

But this is easily checked, since $P_m(u)(\xi)=P_m(u)(R\xi)$ for all $\xi\in\R^2$.
Passing to polar cordinates, we write $P_m(u)(\xi_1,\xi_2)\equiv\tilde P_m(\rho,\theta)$,
with $\xi_1=\rho\cos\theta$ and $\xi_2=\rho\sin\theta$ and observe that  for
each fixed $\rho>0$, the trigonometric polynomial $\tilde P_m(\rho,\theta)$ has degree
smaller or equal than $m+2$ and period $2\pi/n$.
If $m\le n-3$, it follows that $\tilde P_m(\rho,\theta)\equiv \tilde P_m(\rho,0)$ for all~$\theta$,
{\it i.e.\/} $P_m(u)(\xi)$ is radial.
This implies that $P_m(u)(\xi)$ identically vanishes for odd~$m$ and, for even~$m$, that 
$P_m(u)(\xi)$ has the form $c_m(t)(\xi_1^2+\xi_2^2)^{m/2}$ for some constant $c_m(t)$:
in any case, $P_m(u)(\xi)$ is divisible by $\xi_1^2+\xi_2^2\,$ for $m=0,\ldots,n-3$.  
 
If $n\ge3$, then Proposition~\ref{proposition1} applies with $M=n-3$ and
we deduce that, for all $T\ge0$,  $u\in C([0,T],L^\infty_{n+1}(\R^2))$.

\medskip
Solutions that are invariant just under the group $C_n$, in general, decay slowly as $t\to\infty$.
To obtain solutions with large decay rates in space-time we will put more symmetries 
on the data and make use of the vorticity formulation.
Assume now that {\em the flow is invariant
under the dihedral group $D_n$\/} at the beginning of the evolution:
$a(Rx)=Ra(x)$, $a(\tau x)=\tau a(x)$, and that $\omega_0=\partial_1 a_2-\partial_2 a_1$ is rapidly
decreasing as $|x|\to\infty$.
If $\omega(x,t)$ is the unique solution to the two-dimensional vorticity equation starting from~$\omega_0$,
then $\omega(x,t)=\omega(Rx,t)=-\omega(\tau x,t)$, for all $x\in\R^2$ and $t\ge0$ and, by
Proposition~\ref{proposition2}, $\omega(x,t)$ is rapidly decreasing as $|x|\to\infty$ for all fixed $t\ge0$.
It now remains to compute the number of vanishing moments of $\omega(x,t)$:

\begin{lemma}
\label{lemma2}
If $\omega(x,t)$ is as above, then the moments $\int x^\alpha\omega(x,t)\,dx$ vanish
for all $t\ge0$, and all double-index $\alpha=(\alpha_1,\alpha_2)$ such that 
$|\alpha|=\alpha_1+\alpha_2\le n-1$
\end{lemma}

\Proof
Taking the Fourier transform in the space variables, we see that $\widehat\omega(\xi,t)\in C^\infty(\R^2)$
for all $t\ge0$ and $\widehat\omega(\xi,t)=\widehat\omega(R\xi,t)$, $\widehat\omega(\xi,t)=-\widehat\omega(\tau\xi,t)$.
In particular, $\widehat\omega(\xi,t)$ identically vanishes on $n$ different straight lines passing through $\xi=0$.
The Taylor formula then implies that $\partial^\alpha_\xi\widehat\omega(0,t)=0$ for all $\alpha$ such that
$|\alpha|\le n-1$ and Lemma~\ref{lemma2} follows.
\endProof

We can summarize the results of this section in the following theorem, which sharpen
the conclusion of Theorem~\ref{theorem1}.

\begin{theorem}
\label{theorem3}
Let $a=(a_1,a_2)$ be a rapidly decreasing and divercence-free vector field in $\R^2$.
If $a$ is invariant under the cyclic group $C_n$ (n=3,4\ldots)  then the strong solution $u(x,t)$
to (NS) such that $u(0)=a$ satisfies $u(x,t)=O(|x|^{-(n+1)})$ as $|x|\to\infty$ for all $t\ge0$.

If, in addition, $a$ is invariant under the dihedral group $D_n$ and $\omega_0=\partial_1 a_2-\partial_2 a_1$
is also rapidly decreasing in $\R^2$, then the moments of the vorticity $\omega(x,t)$ of the flow
vanish up to the order $n-1$ for all $t\ge0$, 
$\sup_{t\ge0}|u(x,t)|\le C(1+|x|)^{-(n+1)}$
and \eqref{profiles-omega}-\eqref{Lp-decay-omega} and
\eqref{Lp-decay-u-2d} hold true. 
\end{theorem}


\begin{remark}
In the case $n=4$, the symmetries described in the second part of Theorem~\ref{theorem3},
are the same as those studied in \cite{Bra01}: indeed the fact that the flow is invariant under 
the dihedral group $D_4$, can be written as follows:
$u_1(x_1,x_2,t)=-u_1(-x_1,x_2,t)$, $u_1(x_1,x_2,t)=u_1(x_1,-x_2,t)$ and $u_1(x_1,x_2,t)=u_2(x_2,x_1,t)$,
which are exactly the conditions of \cite{Bra01} in the two dimensional case.
\end{remark}

\section{The three dimensional case}
\label{section3D}

We now study the class of flows which are invariant under finite 
subgroups of the group of all the isometries of  the space.
We shall identify two of such groups $G$ and $G'$ if they are conjugate in $O(3)$
({\it i.e.\/} $G\sim G'$ if there exists an orthogonal matrix~$T$ such that $G'=TGT^{-1}$).
Note that two flows that that are invariant under groups which are isomorphic, but not conjugate,
may behave quite differently and this is why will not identify groups which are simply isomorphic
in what follows.

\subsection{Finite subgroups of $O(3)$}

\paragraph{Finite groups of proper rotations}
The material of this section is very classical, but we present it to fix some notations.
We start recalling the well known classification of all finite subgroups of the special orthogonal
group $SO(3)$.
We closely follows the presentation given in \cite{Mil72}.
If ${\cal S}$ is any subset of $\R^3$, the group $G({\cal S})$ of all $P\in SO(3)$ such that $P$ leaves ${\cal S}$
globally invariant is called {\em the complete direct symmetry group of ${\cal S}$\/}.
For different choices of ${\cal S}$ we obtain in this way only five different types of groups that are listed below:
For each group we shall indicate a set of matrices generating
$G({\cal S})$
since we will need these generators in our subsequent calculations. 
\begin{enumerate}
\item
If ${\cal S}$ is a {\em $n$-pyramid\/}, $n=1,2,\ldots$
({\it i.e.\/} a right piramid with base a $n$-sided regular
polygon such that the distance from the vertex of the pyramid to a vertex of the base is not equal to 
one side of the polygon, with an obvious modification if $n=1,2$),
then $G({\cal S})$ is the cyclic group $C_n$ of order~$n$.
A generator of this group is {\it e.g.\/}
\begin{equation}
\label{Rn}
 R_n=\begin{pmatrix}
  \cos(2\pi/n) & -\sin(2\pi/n) & 0\\
  \sin(2\pi/n) & \cos(2\pi/n) & 0\\
  0 & 0 & 1
  \end{pmatrix}.
\end{equation}
\item
If ${\cal S}$ is a {\em $n$-prism\/}, $n=2,3,\ldots$ 
({\it i.e.\/} a right cylinder with base a $n$-sided regular
polygon and height not equal to one side of the polygon, modification if $n=2$),
then $G({\cal S})$ is the dihedral group $D_n$ of order $2n$.
This group is generated by $R_n$ and by
\begin{equation}
\label{U}
 U=\begin{pmatrix}
  1 & 0 & 0\\
  0 & -1 & 0\\
  0 & 0 & -1
  \end{pmatrix}.
\end{equation}
\item
If ${\cal S}$ is a {\em tetrahedron\/}, then $G({\cal S})={\bf T}$ (the tetrahedral group).
This group has order $12$, is isomophic to the alterning group $A_4$ and it is generated by
a rotation by $2\pi/3$ around an axis passing through a vertex and the center of ${\cal S}$
and by a rotation by $\pi$ around an axis passing through the midpoints of two opposite edges.
If $(-1,-1,-1)$, $(1,1,-1)$, $(-1,1,1)$ and $(1,-1,-1)$ are the vertices of ${\cal S}$, then we see
that two generators of ${\bf T}$ are $U$ and
\begin{equation}
S=\begin{pmatrix}
  0 & 0 & 1\\
  1 & 0 & 0\\
  0 & 1 & 0
  \end{pmatrix}.
\end{equation}
\item
If ${\cal S}$ is a {\em cube (or an octahedron)\/} then $G({\cal S})={\bf O}$ (the octahedral group).
This group has order~$24$ and is isomorphic to the symmetric group $S_4$.
If $(\epsilon_1,\epsilon_2,\epsilon_3)$, ($\epsilon_j=1$ or $-1$, $j=1,2,3$) are the vertices of the cube,
then we see that ${\bf O}$ is generated by $U$, $S$ and
a rotation $V$ by $\pi$ around an axis passing through the midpoints of two opposite edges of the cube.
We may choose
\begin{equation}
\label{V}
 V=\begin{pmatrix}
  0 & 1 & 0\\
  1 & 0 & 0\\
  0 & 0 & -1
  \end{pmatrix}. 
\end{equation}
We finally observe that ${\bf T}$ is a subgroup of index~$2$ in ${\bf O}$.

\item
If ${\cal S}$ is an {\em icosahedron (or a dodecahedron)\/} then
$G({\cal S})={\bf Y}$ (the icosahedral group).
This group has order~$60$, is isomorphic to the alterning group $A_5$, it is generated
by a rotation by $2\pi/5$ around an axis passing through two opposite vertices of the icosahedron and 
a rotation by $2\pi/3$ around an axis passing through the center of two opposite faces.

If $(\pm\phi,0,\pm1)$, $(0,\pm1,\pm\phi)$ and $(\pm1,\pm\phi,0)$ are the~$12$ vertices of ${\cal S}$
(here $\phi=(\sqrt5-1)/2$ is the gold number) then we see that ${\bf Y}$ is generated by
$S$ and
\begin{equation}
\label{J}
 J=\begin{pmatrix}
  {1\over2} & -{\sqrt{5}+1\over 4} & {\sqrt{5}-1\over4} \\
  {\sqrt{5}+1\over 4} & {\sqrt{5}-1\over4} & -{1\over2} \\
  {\sqrt{5}-1\over4} & {1\over2} &  {\sqrt{5}+1\over 4}
  \end{pmatrix}.
\end{equation}
It also easily seen that ${\bf Y}$ contains ${\bf T}$. Indeed, the transfomation~$U$
corresponds now to a rotation by~$\pi$ around an axis passing through the midpoints of to opposite edges of
the icosahedron.
\end{enumerate}

\begin{remark}
A classical result states that if $G$ is a finite subgroup of $SO(3)$ then $G$ is
conjugate to one the preceding five groups.
For more details on those groups we refer {\it e.g.\/} to \cite{Mil72}.
\end{remark}

\paragraph{Other finite groups of isometries.}
It ${\cal S}$ is a subset of $\R^3$, then {\em the complete symmetry group of ${\cal S}$\/} si defined as
the group of all orthogonal transformations which leave ${\cal S}$ globally invariant.
Let us recall that if $G$ is a finite subgroup of $O(3)\backslash SO(3)$, such that the inversion $I$
(the symmetry with respect to the origin) belongs to $G$, then $G=G_1\cup IG_1$, where $G_1=G\cap SO(3)$ is
one of the five groups of proper rotation considered in the preceding paragraph.
In this case $G$ is obtained as direct product of $G_1$ and a cyclic group of order~2, and 
the generators of $G$ are the same of $G_1$, together with 
\begin{equation}
\label{I}
 I=\begin{pmatrix}
  -1 & 0 & 0\\
  0 & -1 & 0\\
  0 & 0 & -1
  \end{pmatrix}.
\end{equation}

On the other hand, if $G$ is a finite subgroup of $O(3)\backslash SO(3)$, but $I$ does not belong to $G$,
then $G^+\equiv G_1\cup \{Ig\,:\, g\in G\backslash G_1\}$
is a finite subgroup of $SO(3)$, which contains a finite subgroup $G_1$ of index~$2$.
Further, $G^+$ is isomorphic
(but not conjugate) to $G$. The group $G$ is usually noted $G^+G_1$ in the literature of point groups.
We thus can form four more types of group in this way, namely
$C_{2n}C_n$, $D_nC_n$, $D_{2n}D_n$ and ${\bf O}{\bf T}$.

We now follow the classical classification of Sch\"oenflies (see also \cite{Mil72}, \cite{Yal88}), starting with
the complete symmetry groups of suitably modified prisms.
\begin{enumerate}
\item
we lump togheter the groups $C_n\cup IC_n$
for odd $n$ with the groups $C_{2n}C_n$ for even $n$, to form {\em the cyclic group\/} ${\bf S}_{2n}$ of order $2n$
(the complete symmetry group of an {\em alternating $2q$-prism}, see \cite{Yal88} for a plot).
This group is generated by a rotation-inversion by $\pi/n$
(a rotation of $\pi/n$ around an axis followed by a reflection with respect to a plane perpendicular to the axis):
\begin{equation}
\label{RIn}
 \tilde R_{n/2}=\begin{pmatrix}
  \cos(\pi/n) & -\sin(\pi/n) & 0\\
  \sin(\pi/n) & \cos(\pi/n) & 0\\
  0 & 0 & -1
  \end{pmatrix}.
\end{equation}
\item
Lumping together the groups $C_n\cup I C_n$
for even~$n$ with the groups $C_{2n}C_n$ for odd $n$, we form an {\em abelian group \/}
of order~$2n$, denoted by ${\bf C}_{nh}$
(the complete symmetry group of a {\em shaved $q$-prism\/}).
This group is generated by a rotation-inversion $\tilde R_n$ by~$2\pi/n$
and a rotation $R_n$
by $2\pi/n$  around the same axis.
Note that ${\bf C}_{nh}$ turns out to be a cyclic group if $n$ is odd, but this group is not conjugate to ${\bf S}_{2n}$.
\item
The group $D_nC_n$ has order~$2n$ and is usually denoted by ${\bf C}_{nv}$.
This is the complete symmetry group of a $n$-pyramid and is formed by $n$ rotations by multiples of $2\pi/n$
around the axis of the pyramid and $n$ reflections in~$n$ vertical planes passing through this axis.
A system of generators of ${\bf C}_{nv}$ is  $R_n$ and
\begin{equation}
\label{W2}
 {W_2}=\begin{pmatrix}
  1 & 0 & 0\\
  0 & -1 & 0\\
  0 & 0 & 1
  \end{pmatrix}.
\end{equation}
\item
Combining $D_{2n}D_n$ for odd $n$ with $D_n\cup I D_n$ for even~$n$ forms the complete group
of symmetry of a $n$-prism, which is denoted by ${\bf D}_{nh}$. This group has order~$4n$, and it is generated
by $R_n$, $W_2$ and by
\begin{equation}
\label{W3}
 {W_3}=\begin{pmatrix}
  1 & 0 & 0\\
  0 & 1 & 0\\
  0 & 0 & -1
  \end{pmatrix}.
\end{equation}
\item
Combining $D_{2n}D_n$ for even $n$ with $D_n\cup I D_n$ for odd~$n$ forms the 
group ${\bf D}_{nd}$, which is the complete symmetry group of a {\em twisted $n$-prism\/}
(the solid obtained pasting two $n$-prisms at their basis, in a such way that the prisms
are rotated by $\pi/n$).
This group has order~$4n$ and is generated by $\tilde R_{n/2}$ and $W_2$.
\item
The group ${\bf T}\cup I {\bf T}$ is denoted by ${\bf T}_h$.
This group has order~$24$, is isomorphic to $A_4\times \Z/2\Z$ and is generated by $S$, $U$ and $ I $
(or simply by $S$ and $W_2$).\footnote
{
The ``symmetric solutions'' $u(x,t)$ introduced in \cite{Bra01} are precisely the flows
which are invariant under the group ${\bf T}_h$.
These solutions have been later considered in \cite{FM02}, \cite{GW2}, \cite{Mi02}, \cite{Mi02II} and \cite{Bra02},
but the connection with the group ${\bf T}_h$ does not seem to have been noticed.
}
The group ${\bf T}_h$ corresponds to the complete symmetry group of a solid obtained from a cube
shaving off the eight vertices (this solid is often called {\em modified cube\/}, see \cite{Yal88}).
\item
The group ${\bf O}{\bf T}$ is denoted by ${\bf T}_d$. This is the complete symmetry group of a tetrahedron, it
has order $24$, is isomorphic to ${\bf O}$ (hence to $S_4$), but ${\bf T}_d$ and ${\bf O}$ are not conjugate.
This group is generated by the two generators $S$ and $U$ of ${\bf T}$, together with
a reflection $Z$ with respect to a plane passing through  the midpoint of an edge and containing the opposite edge
of the tetrahedron:
\begin{equation}
\label{Z}
 Z=\begin{pmatrix}
  0 & 1 & 0\\
  1 & 0 & 0\\
  0 & 0 & 1
  \end{pmatrix}.
\end{equation}
\item
The group ${\bf O}\cup I {\bf O}$, denoted by ${\bf O}_h$, is the complete groupe of symmetry of a cube
(and of an octahedron).\footnote
{It was pointed out by Kida~\cite{Kid85I} that it is possible to construct
solutions that are both invariant under ${\bf O}_h$ and $2\pi$-periodic
in any direction.
}
 This group is isomorphic to $S_4\times \Z/2\Z$ and contains the~$48$
orthogonal matrices formed by $0$, $1$ and $-1$.
A system of generators for ${\bf O}_h$ is {\it e.g.\/} $S$, $V$ and $ I $.
\item
The group ${\bf Y}\cup I {\bf Y}$ is denoted by ${\bf Y}_h$, and it
is the complete groupe of symmetry of an icosahedron (and of a dodecahedron).
This group has order~$120$, is isomorphic to $A_5\times \Z/2\Z$ and is generated by
$S$, $J$ and $ I $.
\end{enumerate}

\subsection{Application to the Navier--Stokes equations}
\label{application3D}

\paragraph{Space decay}
This paragraph is devoted to the computation of the space decay rates of flows $u(x,t)$
which are invariant under a discrete subgroup of $O(3)$
and such that $u(x,0)$ is localized.
We will not consider here all the possible groups $G$ listed in the preceding section, but we will just
treat the case in which $G$ is either ${\bf T}$, ${\bf T}_h$, ${\bf O}$ or ${\bf Y}$.
Combining the results of this paragraph with the examples of section~\ref{optimality},
however, will immediately give the optimal space decay rates for {\em all\/} groups.

We need the following lemma.

\begin{lemma}
\label{lemma3}
Let $P_m(\xi)$, where $\xi=(\xi_1,\xi_2,\xi_3)\in\R^3$, be
a homogenous polynomial of degree $m+2$ ($m=0,1,\ldots$).
\begin{enumerate}
\item 
If $P_0$ is invariant under the transformations of the tetrahedral group ${\bf T}$, then 
$P_0(\xi)\equiv c_0(\xi_1^2+\xi_2^2+\xi_3^2)$.
\item 
If $P_1$ is invariant under the transformations of either ${\bf O}$, ${\bf T}_h$ or ${\bf Y}$,
then $P_1(\xi)\equiv0$.
\item
If $P_2$ and $P_3$ are invariant under ${\bf Y}$, then $P_2(\xi)\equiv  c_2(\xi_1^2+\xi_2^2+\xi_3^2)^2$
and $P_3(\xi)\equiv0$.
\end{enumerate}
\end{lemma}

\Proof
The proof follows by imposing $P_m(\xi)\equiv P_m(Q\xi)$ where, in the first case $Q=S,U$;
in the second case we take, respectively,
$Q=S,U,V$,  $Q=S,W_2$, or $Q=S,J$; in the third case we choose
$Q=S,J$.
We thus obtain linear systems where the unknowns are the coefficients of $P_m(\xi)$.
Conclusion of Lemma~\ref{lemma3} then immediately follows from lengthy but elementary calculations.
\endProof

Applying Proposition~\ref{proposition1} we immediately get the following

\begin{corollary}
\label{corollary1}
Let $a=(a_1,a_2,a_3)$ be a soleinoidal and rapidly decreasing vector field in $\R^3$
and $u(x,t)$ the strong solution to (NS), which is defined in some time interval $[0,T]$ ($T>0$) ,
such that $u(0)=a$.
If $a$ is invariant under the transformations of ${\bf T}$,
then $u(x,t)=O(|x|^{-5})$ as $x\to\infty$ uniformely in $t\in[0,T]$.
Such decay rate is improved up to $u(x,t)=O(|x|^{-6})$, if $a$ is invariant under either~${\bf O}$ or
${\bf T}_h$, and up to $u(x,t)=O(|x|^{-8})$ if~$a$ is invariant under~${\bf Y}$.
\end{corollary}

\begin{remark}
Note that the homogeneous polynomial $P_1(\xi)\equiv \xi_1\xi_2\xi_3$ satisfies $P_1(\xi)=P_1(S\xi)=P_1(U\xi)$
for all $\xi$. This polynomial is then invariant under ${\bf T}$, but it is not divisible
by $\xi_1^2+\xi_2^2+\xi_3^2$.
On the other hand, the polynomial $P_2(\xi)\equiv \xi_1^4+\xi_2^4+\xi_3^4$
is invariant under both ${\bf O}$ and ${\bf T}_h$, but it is not divisible by $\xi_1^2+\xi_2^2+\xi_3^2$.
At the same way, it is not difficult to construct a homogeneous polynomial of degree~$6$, which
is invariant under ${\bf Y}$ and is not divisible by  $\xi_1^2+\xi_2^2+\xi_3^2$:
we may take {\it e.g.\/} 
$P_4(\xi)\equiv \xi_1^6+\xi_2^6+\xi_3^6 +{3\over4}(5+\sqrt{5})(\xi_1^4\xi_2^2+\xi_2^4\xi_3^2+\xi_3^2\xi_2^2)
+{3\over4}(5-\sqrt{5})(\xi_1^4\xi_3^2+\xi_2^4\xi_1^2+\xi_3^2\xi_2^2)$.

These considerations show that the decay rates computed in Corollary~\ref{corollary1} seem to be optimal
for generic flows invariant under one of the four preceding groups.
We will see by means of the examples of section~\ref{optimality}
that (a) the decay $|x|^{-5}$ is indeed optimal, in general, for flows which are
invariant under the group ${\bf T}_d$, which contains ${\bf T}$, 
(b) the decay $|x|^{-6}$ is optimal inside the group ${\bf O}_h$ which contains both ${\bf O}$ and ${\bf T}_h$,
(c) the decay $|x|^{-8}$ is optimal inside the group ${\bf Y}_h$ which contains ${\bf Y}$.
Finally, we will see that flows which are invariant under the complete group of symmetries $D_{2nh}$ of a
$2n$-prism (which contains all the other groups $C_n$, $S_{2n}$, $D_n$, $C_{nh}$, $C_{nv}$,
$D_{nd}$ and $D_{nh}$), in general, do not decay  faster  than $|x|^{-4}$.
This provides a complete answer to the space decay problem of flows with this kind of symmetries.
\end{remark}

\paragraph{Time decay}
We now compute the time decay rate of flows invariant
under the complete symmetry group of the solids listed in the preceding section.
We will detailed arguments only for the group of the icosahedron ${\bf Y}_h$,
since this group provides the largest decay rates.
Of course similar (but simpler!) considerations can be repeated for the other groups.

As in the two-dimensional case, we shall make use of the vorticity formulation.
We start observing that requiring the condition $P\,u(x,t)\equiv u(P\,x,t)$, for a given $P\in O(3)$,
is equivalent, at least when the singular integral
\eqref{BS3} makes sense, to requiring that
\begin{equation}
\label{inv-Omega}
P\,\Omega(x,t) =\det(P)\,\Omega(P\,x,t)
\end{equation}
for all $x\in\R^n$ and $t\ge0$.
A simple way to prove the equivalence between \eqref{inv-u} and
\eqref{inv-Omega} is to use the Fourier transform and the identity $(Pv)\times(Pw)=\det(P)P(v\times w)$,
which holds for all $v,w\in\R^3$ and $P\in O(3)$.

From now on we shall assume that $a(x)$ is invariant under the group ${\bf Y}_h$ and that
the initial vorticity $\Omega_0=\nabla\times a$ is a rapidly decreasing  vector field as $|x|\to\infty$.
Then we have the following.

\begin{lemma}
\label{lemma4}
Let $\Omega$ be a rapidly decreasing vector field in $\R^3$, such that $P\Omega(x)\equiv\det(P)\Omega(Px)$
for all transformations $P$ belonging to the complete symmetry group of the icosahedron.
Then the moments of $\Omega$ vanish up to the order~$5$.
\end{lemma}

\Proof
Let us denote by $P_j\in {\bf Y}_h\subset O(3)$ the reflection with respect to a plane $\pi_j$ of symmetry of the
icosahedron  $(j=1,\ldots,15)$.
Note that each of the six axes passing through two opposite vertices of the icosahedron belongs exactly
to five distinct planes. 
Let ${\bf e}_1,\ldots,{\bf e}_6 $ be six unit vectors corresponding to 
these axes and such that ${1\over\sqrt2}({\bf e}_1+{\bf e}_2)=(1,0,0)$, 
${1\over\sqrt2}({\bf e}_3+{\bf e}_4)=(0,1,0)$ and ${1\over\sqrt2}({\bf e}_5+{\bf e}_6)=(0,0,1)$
(this is possible if we choose the vertices of the icosahedron as in the preceding section).
If we show that $\int x^\alpha \langle \Omega(x),{\bf e}_k\rangle=0$ for $k=1,\ldots,6$,
and some $\alpha\in\N^3$ (where $\langle\cdot,\cdot\rangle$ denotes the scalar product in $\R^3$),
then it will follow that $\int x^\alpha\Omega(x)\,dx=0$.

Condition $P_j\Omega(x)\equiv -\Omega(P_j x)$ implies that $\Omega(x)$ is orthogonal to $\pi_j$,
for all $x\in \pi_j$.
In the same way, passing to the Fourier transform we see that 
$\widehat\Omega(\xi)$ is orthogonal to $\pi_j$ for all $\xi\in\pi_j$.
In particular, for each $k=1,\ldots,6$
there exist five planes containing the axis generated by ${\bf e}_k$, on which the function
$f_k(\xi)\equiv\langle\widehat\Omega(\xi),{\bf e}_k\rangle$
identically vanish.

Now we use the general fact that if $g(\xi)\in C^\infty(\R^3)$ identically vanishes on~$n$ distinct planes
passing through a given axis, then $g$ has vanishing derivatives on this axis up to the order~$n-1$
(this simple fact can be seen using the same argument as in Lemma~\ref{lemma2}).
 
Since $f_k$ is a smooth function, it follows that the derivatives of $f_k$ identically vanish
up to the order~$4$ on the $k$-th axis.
This shows that the moments of $\Omega$ vanish up to the order~$4$.

But the group ${\bf Y}_h$ contains the three reflections with respect to the planes
$x_1=0$, $x_2=0$ and $x_3=0$ (assuming that the icosahedron is orientated as above).
Therefore $\Omega_i(x_1,x_2,x_3)$ is an even function with respect to $x_i$ and an odd function
with respect to $x_h$ ($i,h=1,2,3$ and $i\not=h$).
It then follows that for any $\alpha\in\N^3$, such that $\alpha_1+\alpha_2+\alpha_3$ is an odd
integer, $\int x^\alpha\Omega_i(x)\,dx=0$ $(i=1,2,3)$.
Lemma~\ref{lemma4} is thus proved.
\endProof

Combining this result with Proposition~\ref{proposition1} and Lemma~\ref{lemma1} implies
the following:

\begin{corollary}
\label{corollary2}
Let $\Omega_0$ be a divergence-free and rapidly decreasing vector field,
such that $P\Omega_0(x)=\det(P)\Omega_0(Px)$ for all $P\in {\bf Y}_h$.
If $\sup_x |x|^2|\Omega_0(x)|$ is small, then the solution
$\Omega(x,t)$ of Proposition~\ref{proposition2} satisfies
\eqref{profiles-Omega}-\eqref{Lp-decay-Omega} with $n=6$.
Furthermore, the corresponding velocity field belongs to $C([0,+\infty[,L^\infty_8(\R^3))$
and satisfies \eqref{Lp-decay-u-3d} (with $n=6$).
\end{corollary}

\medskip
These arguments apply also to the simpler case of flows invariant to complete symmetry group
of the tetrahedron and the complete symmetry group of the cube
(and, with slight modification, to their subgroups).
This yields {\it e.g.\/} that \eqref{Lp-decay-u-3d} holds true with $n=3$ in the case of flows invariant
under ${\bf T}_d$, and with $n=4$ in the case of flows invariant
under ${\bf O}_h$.
We leave the corresponding computations to the reader.
Theorem~\ref{theorem2} then follows.

\section{Examples of localized flows}
\label{optimality}

In this section we provide explicit examples of initial data leading to
flows invariant under the groups considered in the preceding section.
These examples also shows that the space decay rates previously computed are sharp.
The proof of the optimality is based on the following fact: if $a(x)$ is a rapidly decreasing
divergence-free vector field in $\R^d$ $(d\ge2)$, such that the homogeneous polynomial
$P_m(a)(\xi)$, $\xi=(\xi_1,\ldots,\xi_d)\in\R^d$, defined by~\eqref{Pma} {\em is not divisible\/}
by $\xi_1^2+\ldots+\xi_d^2$, then there exists a decreasing sequence $t_k\to0$ such that
\begin{equation}
\label{m-spreading}
\limsup_{R\to\infty} R^{1+m}\int_{R\le|x|\le 2R} |u(x,t_k)|\,dx >0 \qquad\mbox{for all $k=1,2,\ldots$},
\end{equation}
where $u(x,t)$ is the strong (local) solution to (NS) starting from $a$
(see \cite{BrM}).
Condition~\eqref{m-spreading} implies that $u$ cannot decay faster than $|x|^{-(d+1+m)}$ 
uniformely in any time interval $[0,T]$ $(T>0)$.

From now on we shall assume $d=3$.
To give examples of soleinodal vector fields $a$ which are invariant under an orthogonal transformation
$P$, it will be often convenient first to construct a {\em potential vector\/} field $b$, such that
$b(Px)=\det(P)Pb(x)$ and then to define $a=\nabla \times b$.

Let $G$ be a finite subgroup of $O(3)$ which is not of polihedral type ({\it i.e.\/} $G$ does not
contain ${\bf T}$).
If $G$ has order~$n$ then it is contained in either $D_{nd}$ or in $D_{nh}$.
These two groups are in turn both contained in the complete symmetry group of a $2n$-prism
$D_{2nh}$, which is of order~$8n$.
Let us show that generic flows which are invariant under $D_{2nh}$ do not decay faster than $|x|^{-4}$.
This is immediate: we can take {\it e.g.\/} a vector field of the form
\begin{equation}
a(x)=(-\partial_2\mu(x),\partial_1\mu(x),0),
\end{equation}
where $\mu\in{\cal S}(\R^3)$ is a non-trivial function such that 
$\mu(x_1,x_2,x_3)=-\mu(x_1,-x_2,x_3)=-\mu(x_1,x_2,-x_3)$ and $\mu$ is invariant under a rotation of $\pi/n$
around the vertical axis.
Then $a$ is invariant under $D_{2nh}$, but $\int a_1^2(x)\,dx\not=\int a_3(x)^2\,dx$, hence
\eqref{m-spreading} holds true with $m=0$.
A slight modification of the choice of $\mu$ would show, in the same way, that flows invariant
under the complete group of direct symmetry of the cylinder do not decay faster than $|x|^{-4}$,
in general.

A very simple example of a vector field which is invariant under the complete symmetry group ${\bf T}_d$ of
a tetrahedron is obtained choosing {\it e.g.\/} 
the potential vector $\bar b_1(x_1,x_2,x_3)=x_1(x_2^2-x_3^2)e^{-|x|^2}$,
$\bar b_2(x_1,x_2,x_3)=\bar b_1(x_2,x_3,x_1)$ and $\bar b_3(x_1,x_2,x_3)=\bar b_1(x_3,x_1,x_2)$.
Another possible simple choice for the first component of the potential vector would be 
$\tilde b_1(x_1,x_2,x_3)=x_2x_3(x_2^2-x_3^2)e^{-|x|^2}$
These two choices give, respectively,
\begin{equation}
\bar a(x)=\begin{pmatrix} -2x_2x_3(2+2x_1^2-x_2^2-x_3^2)e^{-|x|^2}\\
	-2x_3x_1(2+2x_2^2-x_3^2-x_1^2)e^{-|x|^2}\\
	-2x_1x_2(2+2x_3^2-x_1^2-x_2^2)e^{-|x|^2} 
 \end{pmatrix}
\end{equation}
and
\begin{equation}
\tilde a(x)=\begin{pmatrix}
	x_1(2x_1^2- 3x_2^2- 2x_1^2x_2^2+ 2x_2^4- 3x_3^2+ 2x_3^4- 2x_1^2x_3^2) e^{-|x|^2} \\
	x_2(2x_2^2- 3x_3^2- 2x_2^2x_3^2+ 2x_3^4- 3x_1^2+ 2x_1^4- 2x_2^2x_1^2) e^{-|x|^2} \\
	x_3(2x_3^2- 3x_1^2- 2x_3^2x_1^2+ 2x_1^4- 3x_2^2+ 2x_2^4- 2x_3^2x_2^2) e^{-|x|^2}
\end{pmatrix}
\end{equation}
Accordingly with Lemma~\ref{lemma2}, the two polynomials $P_0(\bar a)(\xi)$ and $P_0(\tilde a)(\xi)$
are divisible by $\xi_1^2+\xi_2^2+\xi_3^3$.
Hence the solutions $\bar u$ and $\tilde u$ starting from $\bar a$ and $\tilde a$
decay at least as fast as $|x|^{-5}$, at the beginnning of their evolution.
However, a direct calculation shows that both $P_1(\bar a)(\xi)$ and $P_1(\tilde a)(\xi)$ identically vanish.
This means that $\bar u$ and $\tilde u$ may decay faster than expected for generic flows invariant under
the group~${\bf T}_d$.
However, one easily checks that $P_1(\bar a+\tilde a)(\xi)\equiv c\xi_1\xi_2\xi_3$ for some
constant $c\not=0$.
Hence, the flow starting from $(\bar a+\tilde a)(x)$, which is also invariant under ${\bf T}_d$
(and, in particular, under ${\bf T}$), cannot decay faster than $|x|^{-5}$.
This decay rate is thus sharp, in general, for the groups ${\bf T}$ and ${\bf T}_d$.

\smallskip
Note that the field $\tilde a(x)$ turns out to be invariant under the transformation of the larger
group ${\bf O}_h$.
Now, it is not difficult to check that the homogeneous polynomial 
$P_2(\tilde a)(\xi)$ is not divisible by $\xi_1^2+\xi_2^2+\xi_3^2\,$.
Indeed, a necessary condition on the coefficients of $P_2(\tilde a)(\xi)$, to obtain a polynomial
which is divisible by $\xi_1^2+\xi_2^2+\xi_3^2$, would be 
$4\int x_1x_2\tilde a_1\tilde a_2(x)\,dx=
 \int(x_1^2-x_2^2)(\tilde a_1^2-\tilde a_2^2)(x)\,dx$.
But the left hand side equals ${57\over 512}\pi^{3/2}\sqrt2$ and the right-hand side equals
${15\over 64}\pi^{3/2}\sqrt2$.
Then, the decay rate $|x|^{-6}$ is optimal for generic flows invariant under
the complete symmetry group of a cube
(hence, also for the groups ${\bf O}$ and ${\bf T}_h$ which are both contained in ${\bf O}_h$).

\smallskip
We do not give explicit examples of localized and soleinoidal vector fields invariant under the complete
symmetry group of the icosahedron, since they have a quite complicated expression.
Their existence, however, is obvious: indeed, the fiveteen planes of symmetry of the icosahedron
divide $\R^3$ into~$120$ congruent pyramidal regions (each of them is the convex hull of three
half straight lines arising from the origin).
If $\Gamma$ is one of these regions, it is then sufficient to construct a vector field $a$
which is localized and divergence-free in~$\Gamma$ and such that $\langle a(x),{\bf n}\rangle=0$
on the boundary of $\Gamma$ (${\bf n}$ denotes here the exterior normal).
The extension of such field by subsequent
reflections with respect to the fiveteen planes is then invariant under the transformations of ${\bf Y}_h$.
Note that the condition on the boundary of $\Gamma$ is conserved by the Navier--Stokes evolution:
this means that the fluid particles will remain in the same region~$\Gamma$ for all time.

It now remains to show that the decay $|x|^{-8}$ is optimal for flows invariant under the group ${\bf Y}_h$
(hence, also for ${\bf Y}$).
We shall only outline the proof: one first observes
that the set  of all homogeneous polynomials of degree~$6$ which are invariant
under ${\bf Y}_h$ form a linear space $H_6$ of dimension~$2$ and 
%
%
%
that the polynomials of $H_6$ which are divisible by $\xi_1^2+\xi_2^2+\xi_3^2$
form a subspace of dimension~$1$.
But it is not difficult to construct (as indicated above)
two different localized initial data $a(x)$ and $a^*(x)$, both invariant under ${\bf Y}_h$
and such that the two polynomials $P_4(a)(\xi)$ and $P_4(a^*)(\xi)$ are linearily independent.
It then follows that at least one of them is not divisible by $\xi_1^2+\xi_2^2+\xi_3^2$.
Condition~\eqref{m-spreading} is then satisfied with $m=4$ and our claim is thus proven.

\bigskip
\paragraph{Acknowledgements}

The author would like to express his gratitude to Y.~Colin de Verdi\`ere and Th.~Gallay
for suggesting him to study the problem treated in this paper
and for their many useful remarks.
He is especially indebted to the former, for
suggesting the statement of Lemma~\ref{lemma3}
which was the starting point of this work.


\end{document}